\documentclass[12pt]{article}

\usepackage{amsmath}
\usepackage{amssymb}
\usepackage{amsfonts}
\usepackage{graphics}

\newtheorem{thmnonumber}{Theorem}
\newtheorem{thm}{Theorem}[section]

\newtheorem{cor}[thm]{Corollary}

\newtheorem{prop}[thm]{Proposition}
\newtheorem{conjecture}{Conjecture}

\newcommand{\be}{\begin{equation}}
\newcommand{\ee}{\end{equation}}
\newcommand{\openbox}{\leavevmode
  \hbox to8pt{\hfil\vrule\vbox to6pt{\hrule width6pt\vfil\hrule}\vrule}}

\newcommand{\qed}{\hbox to5pt{ } \hfill \openbox\bigskip\medskip}

\newcommand{\rk}{\mbox{\rm rank}}

\newcommand{\cF}{\mbox{$\cal F$}}
\newcommand{\cG}{\mbox{$\cal G$}}
\newcommand{\cH}{\mbox{$\cal H$}}

\newcommand{\cM}{\mbox{$\cal M$}}

\newcommand{\ve}[1]{\mathbf{#1}}
\newcommand{\cD}{\mbox{$\cal D$}}

\newcommand{\Z}{\mathbb Z}

\newcommand{\F}{\mathbb F}
\newcommand{\Fe}{\mathbb F _2}

\newcommand{\sh}{\mathrm{sh}}
\newcommand{\Vdim}{\mathop\textup{Vdim}}

\title{An uniform version of Dvir and Moran's theorem}
\author{G\'abor Heged\"{u}s
\\{\normalsize  \'Obuda University}
\\{\normalsize B\'ecsi \'ut 96/B, Budapest, Hungary, H-1032}
\\{\normalsize hegedus.gabor@uni-obuda.hu}
}

\begin{document}

\maketitle
\begin{abstract}
Dvir and Moran proved the following upper bound  for the size of a family $\mbox{$\cal F$}$ of subsets of $[n]$ with $\mathop\textup{Vdim}(\mbox{$\cal F$} \Delta \mbox{$\cal F$})\leq d$.

Let $d\leq n$ be integers. Let $\mbox{$\cal F$}$ be a family of subsets of $[n]$ with $\mathop\textup{Vdim}(\mbox{$\cal F$} \Delta \mbox{$\cal F$})\leq d$. Then
\[
\left|\mbox{$\cal F$}\right|\le 2\sum_{k=0}^{\lfloor d/2 \rfloor}\binom nk.
\]

Our main result is the following uniform version of Dvir and Moran's result.

Let $d\leq n$ be integers. Let $\mbox{$\cal F$}$ be an uniform  family of subsets of $[n]$ with $\mbox{Vdim}(\mbox{$\cal F$} \Delta \mbox{$\cal F$})\leq d$. Then
\[
\left|\mbox{$\cal F$}\right|\le 2  {n \choose \lfloor d/2 \rfloor}.
\]

Denote by $\mathbf{v_F}\in \{0,1\}^n$ the
characteristic vector of a set
$F \subseteq [n]$.

Our proof is based on the following uniform version of Croot-Lev-Pach Lemma:
 
Let $0\leq d\leq n$ be integers. Let $\mbox{$\cal H$}$ be a $k$-uniform  family of subsets of $[n]$.  Let $\mathbb F$ be a field. Suppose that there exists a polynomial $P(x_1, \ldots ,x_n,y_1, \ldots ,y_n)\in \mathbb F[x_1, \ldots ,x_n,y_1, \ldots ,y_n]$ with $\mbox{deg}(P)\leq d$ such that $P(\mathbf{v_F},\mathbf{v_F})\ne 0$ for each $F\in \mbox{$\cal H$}$ and $P(\mathbf{v_F},\mathbf{v_G})= 0$ for each $F\ne G\in \mbox{$\cal H$}$. Then 
\[
\left|\mbox{$\cal H$}\right|\le 2  {n \choose \lfloor d/2 \rfloor}.
\]
\end{abstract}
\medskip

\section{Introduction}

Throughout this paper $n$ denotes a positive integer, and $[n]$ stands for 
the set $\{1,2,\dots, n\}$. We denote 
by $2^{[n]}$ the set of all subsets of $[n]$. Subsets of $2^{[n]}$ are called \emph{set families}.
Let $\binom{[n]}m$ denote the family of all subsets of $[n]$
which have cardinality $m$, and $\binom{[n]}{\le m}$ of all subsets that 
have size at most $m$. 

A family $\cF$ of subsets of $[n]$ is {\em $k$-uniform}, if $|F|=k$ for each $F\in \cF$.

Let $\F$ be a field.
$\F[x_1, \ldots, x_n]=\F[\ve x]$ denotes  the
ring of polynomials in the variables $x_1, \ldots, x_n$ over $\F$.
For a subset $F \subseteq [n]$ we write
$\ve x_F = \prod_{j \in F} x_j$.
In particular, $\ve x_{\emptyset}= 1$.

We denote by $\ve v_F\in \{0,1\}^n$ the
characteristic vector of a set
$F \subseteq [n]$.


It is a challenging old problem to find strong upper bounds for the size of progression-free subsets in finite Abelian groups.
Croot, Lev and Pach achieved recently a breakthrough in this research area and proved a new exponential upper bound for the size of  three-term progression-free subsets  in the groups $({\Z}_4)^n$ (see  \cite{CLP}), where $n\geq 1$ is an arbitrary integer. They based their proof on the following 
simple statement (see \cite{CLP}  Lemma 1).
\begin{prop} \label{CLP}
Suppose that $n \geq 1$ and $d \geq 0$ are integers, $P$ is a multilinear polynomial in $n$ variables of total degree at most $d$ over a field $\F$, and $A \subseteq {\F}^n$ is a subset with 
$$
|A|>2\sum_{i=0}^{d/2} {n\choose i}.
$$
If $P(\ve a- \ve b)=0$ for all $\ve a, \ve b\in A$, $\ve a\neq \ve b$, then $P(\ve 0)=0$.
\end{prop}

Consider a family $\cF$ of subsets of $[n]$. We say that $\cF$ \emph{shatters} 
$M\subseteq[n]$ if
\[
\{F\cap M\,:\,\,F\in \cF\}=2^M.
\]
Define
\[
\sh(\cF)=\{M\subseteq [n]\,:\,\,\cF\text{ shatters }M\}.
\]

\medskip

We say that a family $\cF$ has {\em $VC$-dimension} $m$, if $m$ is the maximum of the size 
of sets shattered by $\cF$. We denote by $\Vdim(\cF)$ the $VC$-dimension of a family $\cF$.

The following result is  fundamental in the theory of shattering.

\begin{thmnonumber}[Sauer{\cite{Sau}}, Perles, Shelah{\cite{Sh}},
Vapnik, Chervonenkis{\cite{VC}}] \label{VCdim}
Let $\cF$ be a family of subsets of $[n]$ with $\Vdim(\cF)\le d$.  Then
\[
\left|\cF\right|\le \sum_{k=0}^{d}\binom nk
\]
and the upper bound is sharp.
\end{thmnonumber}

Let $\cF$ and $\cG$ be  families of subsets of $[n]$. We denote by $\cF \Delta \cG$ the symmetric difference of these families:
$$
\cF \Delta \cG:=\{A\Delta B:~ A\in \cF,B\in \cG\}.
$$ 

Dvir and Moran proved an  upper bound  for the size of a family $\cF$ of subsets of $[n]$ with $\Vdim(\cF \Delta \cF)\leq d$.  Their proof based on Proposition \ref{CLP}.

\begin{thm}\label{dvir}
Let $0\leq d\leq n$ be integers. Let $\cF$ be a family of subsets of $[n]$ with $\Vdim(\cF \Delta \cF)\leq d$. Then
\[
\left|\cF\right|\le 2\sum_{k=0}^{\lfloor d/2 \rfloor}\binom nk.
\]
\end{thm}

Cambie, Gir\~ao  and Kang proved the following improved version of Theorem \ref{dvir} in \cite{CGK}.

\begin{thm}\label{kang}
Let $d<n$ be positive integers with $d\equiv r \pmod 2$ for some 
$r\in \{0,1\}$. Let $\cF$ be a family of subsets of $[n]$ with $\Vdim(\cF \Delta \cF)\leq d$. Then
$$
|\cF|\leq 2^r \sum_{k=0}^{\lfloor d/2 \rfloor} {n-r \choose k}.
$$

\end{thm}
Kleitman's theorem (see \cite{K}) is an immediate consequence of Theorem \ref{kang}. 

\begin{cor} \label{kleitman}
Let $d<n$ be positive integers with $d\equiv r \pmod 2$ for some 
$r\in \{0,1\}$. Let $\cF$ be a family of subsets of $[n]$ with $\cF \Delta \cF\subseteq {[n]\choose {\leq d}}$. Then
$$
|\cF|\leq 2^r \sum_{k=0}^{\lfloor d/2 \rfloor} {n-r \choose k}.
$$
\end{cor}

Our main result is the following uniform version of Theorem \ref{dvir}.

\begin{thm}\label{main}
Let $d\leq n$ be integers. Let $\cF$ be an uniform  family of subsets of $[n]$ with $\Vdim(\cF \Delta \cF)\leq d$. Then
\[
\left|\cF\right|\le 2  {n \choose \lfloor d/2 \rfloor}.
\]
\end{thm}

We prove here a new uniform version of Proposition \ref{CLP}. The proof of Theorem \ref{main} is based completely on this result.

\begin{thm}\label{uniCLP}
Let $0\leq d\leq n$ be integers. Let $\cH$ be a $k$-uniform  family of subsets of $[n]$.  Let $\F$ be a field. Suppose that there exists a polynomial $P(\ve x,\ve y)\in \F[\ve x,\ve y]$ with $\mbox{deg}(P)\leq d$ (here $\ve x=(x_1, \ldots ,x_n), \ve y=(y_1, \ldots ,y_n)$) such that $P(\ve v_F,\ve v_F)\ne 0$ for each $F\in \cH$ and $P(\ve v_F,\ve v_G)= 0$ for each $F\ne G\in \cH$. Then 
\[
\left|\cH\right|\le 2  {n \choose \lfloor d/2 \rfloor}.
\]
\end{thm}

In Section 2 we collected the preliminaries about Gr\"obner  basis theory and standard monomials. In Section 3 we present our proofs. In Section 4 we give an interesting conjecture which strengthens our main result.

\section{Preliminaries}



Define
$V(\cF)$ as the subset  $\{\ve v_F : F \in \cF\} \subseteq \{0,1\}^n \subseteq \F^n$ for any family of subsets $\cF \subseteq 2^{[n]}$.

It is natural to consider the ideal $I(V(\cF))$:
$$ 
I(V(\cF)):=\{f\in \F[\ve x]:~f(\ve v)=0 \mbox{ whenever } \ve v\in V(\cF)\}. 
$$
It is easy to verify that we can identify the algebra $\F[\ve x]/I(V(\cF))$ and the algebra of
$\F$ valued functions on $V(\cF)$. Consequently 
$$
\dim _\F
\F[\ve x]/I(V(\cF))=|\cF|.
$$


We recall some basic facts about  Gr\"obner 
basis theory  and standard monomials.  We refer to \cite{AL}, \cite{CLS} for details.

We say that a linear order $\prec$ on the monomials 
is a {\em term order},  if 1 is
the minimal element of $\prec$, and $\ve u \ve w\prec \ve v\ve w$ holds for any monomials
$\ve u,\ve v,\ve w$ with $\ve u\prec \ve v$. The two most important term orders are the lexicographic
order $\prec_l$ and the deglex order $\prec _{d}$. 
Recall the definition of  the deglex order: we have $\ve u\prec_{d} \ve v$ iff  either
$\deg \ve u <\deg \ve v$, or $\deg \ve  u =\deg \ve v$, and $\ve u\prec_l \ve v$.

The {\em leading monomial} ${\rm lm}(f)$
of a nonzero polynomial $f\in \F[\ve x]$ is the $\prec$-largest
monomial which appears with nonzero coefficient in $f$.

Let $I$ be an ideal of $\F[\ve x]$. A finite subset $\cG\subseteq I$ is a {\it
Gr\"obner basis} of $I$ if for every $f\in I$ there exists a $g\in \cG$ such
that ${\rm lm}(g)$ divides ${\rm lm}(f)$. 
It can be shown that $\cG$ is
actually a basis of $I$, i.e. $\cG$ generates $I$ as an ideal of $\F[\ve x]$ (cf. \cite{CLS} Corollary 2.5.6). 
A well--known fact is (cf. 
\cite[Corollary 1.6.5, Theorem 1.9.1]{AL}) that every
nonzero ideal $I$ of $\F[\ve x]$ has a Gr\"obner basis.

A monomial $\ve z\in \F[\ve x]$ is a {\it standard monomial for $I$} if
it is not a leading monomial for any $f\in I$. We denote by  ${\rm sm}(I)$ 
the set of standard monomials of $I$.

Let $\cF\subseteq 2^{[n]}$ be a set family. 
It is easy to check that the standard monomials of the ideal
$I(\cF):=I(V(\cF))$ are square-free monomials.


It is a fundamental fact that
${\rm sm}(I)$ gives a basis of the $\F$-vector-space $\F[\ve x]/I$.
This means that  
every polynomial $g\in \F[\ve x]$ can be 
uniquely  written in the form  $h+f$ where $f\in I$ and
$h$ is a unique $\F$-linear combination of monomials from ${\rm sm}(I)$. Consequently if $g\in \F[\ve x]$ is an arbitrary polynomial and $\cG$ is a Gr\"obner basis of $I$, then we can reduce $g$ with $\cG$ into a linear combination of
standard monomials for $I$. 

\section{Proofs}

Let $0\leq k\leq n/2$, where $k$ and $n$ are integers. Let $\cM_{k,n}$ stand for the set of all
monomials $\ve x_G$ such that $G=\{s_1<s_2<\ldots <s_j\}\subset [n]$
for which $j\leq k$ and $s_i\geq 2i$ holds for every $i$, $1\leq i\leq j$. 
We write $\cM_k$ instead of the more precise $\cM_{k,n}$, if $n$ is clear from the context. 
It is easy to check that 
$$ 
|\cM_k|={n \choose k}.
$$

Let $\cD_{k,n}$ denote the set of all sets $H=\{s_1<s_2<\ldots <s_j\}\subset [n]$ for which $j\leq k$ and $s_i\geq 2i$ holds for every $i$, $1\leq i\leq j$. 

We described completely the standard monomials  of the complete uniform families of all $k$ element subsets of $[n]$ in \cite{HR}.

\begin{thm}\label{sm1}
Let $\prec$ an arbitrary term order such that $x_n\prec \ldots \prec x_1$. 
Let $0\leq k\leq n$ be integers and define 
$j:=min(k,n-k)$. Then 
$$                                
{\rm sm}(V{[n]\choose k})=\cM_{j,n}.
$$
\end{thm}
\begin{cor} \label{sm3}
Let $0\leq k\leq n$ be integers and define  
$j:=min(k,n-k)$. Suppose that  $d\leq j$. Then 
$$                                
\cD_{k,n}\cap {[n]\choose {\leq d}}=\cD_{d,n}.
$$
\end{cor}

Let $0\leq k\leq n$  be arbitrary integers. Define the vector system
$$
\cF(n,k,2):=V({[n]\choose k})\times V({[n]\choose k})\subseteq \{0,1\}^{2n}.
$$

It is easy to verify the following Corollary from Theorem \ref{sm1}.

\begin{cor} \label{sm2}
Let $\prec$ an arbitrary term order such that  $x_n\prec \ldots \prec x_1$. 
Let $0\leq k\leq n$ be integers and define $j:=min(k,n-k)$. Then 
$$                                
{\rm sm}(\cF(n,k,2))=\{x_{M_1}\cdot y_{M_2}:~ M_1, M_2\in\cD_{j,n} \}\subseteq {\F}[\ve x,\ve y].
$$
\end{cor}

M\'esz\'aros and R\'onyai proved the following result in \cite{MR} Lemma 1 (see also \cite{MR2} Theorem 7).

\begin{thm}\label{smsh}
Let $\prec$ an arbitrary term order such that  $x_n\prec \ldots \prec x_1$.
 Let $\cF$ be a family of subsets of $[n]$. Then ${\rm sm}(V(\cF))\subseteq \{x_U:~ U\in \sh(\cF)\}$.
\end{thm}

{\bf Proof of Theorem \ref{uniCLP}:}

Consider the matrix  $M\in {\F}^{\cH\times \cH}$, where $M_{(F,G)}:=P(\ve v_F,\ve v_G)$ for each $F,G\in \cH$. 

It follows from  the assumptions that  $M$ is a diagonal matrix, where nonzero elements stand in the diagonal, hence 
$$
\rk(M)= |\cH|.
$$

Let $Q$ denote the reduction of $P$ via the deglex Gr\"obner basis of $I(\cF(n,k,2))$. Then $\mbox{deg}(Q)\leq \mbox{deg}(P)\leq d$ and  $M_{(F,G)}=Q(\ve v_F,\ve v_G)$ for each $F,G\in \cH$.  

Let 
$j:=min(k,n-k)$. It follows from Corollary \ref{sm2} that we can write the polynomial $Q$ into the form
$$
Q(\ve x,\ve y)=\sum_{M_1,M_2\in \cD_{j,n} } c_{M_1,M_2} x_{M_1}\cdot y_{M_2}\in \Fe[\ve x,\ve y],
$$
where $c_{M_1,M_2} \in\Fe$ for each $M_1,M_2\in \cD_{j,n}$. After grouping  the terms of the polynomial $Q(\ve x,\ve y)$ we get that
$$
Q(\ve x,\ve y)=\sum_{M\in \cD_{j,n}\cap \{U:~ |U|\leq \lfloor d/2\rfloor\}} c_M \ve x_M g_M(\ve y)+\sum_{J\in \cD_{j,n}\cap \{U:~ |U|\leq \lfloor d/2\rfloor\}} d_J \ve y_J h_J(\ve x),
$$
where $c_M,d_J\in \Fe$, $h_J(\ve x)\in \Fe[\ve x]$, $g_M(\ve y)\in \Fe[\ve y]$ for each $J,M\in \cD_{j,n}$. 

Then it follows from Corollary  \ref{sm3} that
$$
Q(\ve x,\ve y)=\sum_{M\in \cD_{\lfloor d/2\rfloor,n}} c_M \ve x_M g_M(\ve y)+\sum_{J\in \cD_{\lfloor d/2\rfloor,n}} d_J \ve y_J h_J(\ve x),
$$

Since
$$
|\cD_{\lfloor d/2\rfloor,n}|={n \choose \lfloor d/2 \rfloor},
$$
hence we get that
$$
\rk(M)\leq  2  {n \choose \lfloor d/2 \rfloor}.
$$

It follows from the equality $\rk(M)= |\cH|$ that
\[
\left|\cH\right|\le 2  {n \choose \lfloor d/2 \rfloor}.
\]
\qed

{\bf Proof of Theorem \ref{main}:}

Let $\F:=GF(2)$. It follows from Theorem \ref{smsh} that
$$
{\rm sm}(V(\cF \Delta \cF),\prec _{d})\subseteq \{x_U:~ U\in \sh(\cF \Delta \cF)\}.
$$
Since $\Vdim(\cF \Delta \cF)\leq d$, hence 
$$
{\rm sm}(V(\cF \Delta \cF),\prec _{d})\subseteq \{x_U:~ |U|\leq d\}.
$$
Let $\cG$ denote a fixed deglex Gr\"obner basis of $I(V(\cF \Delta \cF))$. Denote by $g:V(\cF \Delta \cF)\to \F$ the function where $g(\ve 0)=1$ and $g(\ve v_T)=0$ for each $T\in \cF \Delta \cF\setminus \{\emptyset\}$. 

If we reduce $g$ with the  Gr\"obner basis $\cG$, we get the polynomial $g'\in \F[\ve x]$. Clearly $\mbox{deg}(g')\leq d$, because $g'$ is a linear combination of deglex standard monomials of $I(V(\cF \Delta \cF))$ and 
${\rm sm}(V(\cF \Delta \cF),\prec _{d})\subseteq \{x_U:~ |U|\leq d\}$.  Since $\cG$ is a Gr\"obner basis of $I(V(\cF \Delta \cF))$, hence 
$$
g(\ve v_G)=g'(\ve v_G)
$$
for each $G\in \cF \Delta \cF$.

Define the polynomial function $f:V(\cF)\times V(\cF)\to \F$ by 
$$
f(\ve x,\ve y):=g'(\ve x+\ve y).
$$
Then
$$
f(\ve v_F,\ve v_F)=g'(\ve 0)=g(\ve 0)=1
$$
for each $F\in \cF$ and 
$$
f(\ve v_F,\ve v_G)=g'(\ve v_F+\ve v_G)=g'(\ve v_{F\Delta G})=g(\ve v_{F\Delta G})=0
$$
for each $F,G\in \cF$, where $F\ne G$.

We can apply Theorem \ref{uniCLP} with the choices $\cH:=\cF$ and $P(\ve x,\ve y):=f(\ve x,\ve y)$. \qed

\section{Concluding remarks}

We think that the next conjecture is the best form of Theorem \ref{main}. 
\begin{conjecture} \label{Hconj}
Let $d<n$ be positive integers with $d\equiv r \pmod 2$ for some 
$r\in \{0,1\}$. Let $\cF$ be an uniform  family of subsets of $[n]$ with  with $\Vdim(\cF \Delta \cF)\leq d$. Then
$$
|\cF|\leq 2^r  {n-r \choose \lfloor d/2 \rfloor}.
$$
\end{conjecture}


\end{document}